%% file: BNR-chap_HAL.tex
\documentclass{imsart}
\usepackage{fixltx2e,fix-cm}
\usepackage{amssymb}
\usepackage{amsmath}
\usepackage{graphicx}
\usepackage{subfigure}
\usepackage{makeidx}
\usepackage{multicol}

\usepackage{color}
\usepackage{hyperref}

\newcounter{chapter}
\newcommand{\chapter}[1]{}

\input{BNR-chap_preamble} 

\begin{document}
\title{On agnostic post hoc approaches to false positive control}

\begin{aug}
\author{\fnms{Gilles} \snm{Blanchard}\ead[label=e1]{gilles.blanchard@math.u-psud.fr}}

\address{
  Universit\"at Potsdam, Institut f\"ur Mathematik\\
  Karl-Liebknecht-Stra{\ss}e 24-25 14476 Potsdam, Germany\\
}

\address{
  Laboratoire de Mathématiques d'Orsay,\\
  Universit\'e Paris-Sud, CNRS, \\
  Universit\'e Paris-Saclay, 91405 Orsay Cedex, France\\
  \printead{e1}
}

\author{\fnms{Pierre} \snm{Neuvial}\ead[label=e2]{pierre.neuvial@math.univ-toulouse.fr}}

\address{
  Institut de Math\'ematiques de Toulouse; \\
  UMR 5219, Universit\'e de Toulouse, CNRS\\
  UPS IMT, F-31062 Toulouse Cedex 9, France\\
  \printead{e2}
}

\author{\fnms{Etienne} \snm{Roquain}\ead[label=e3]{etienne.roquain@upmc.fr}}

\address{
  Sorbonne Universit\'e, Laboratoire de Probabilit\'es, Statistique et Mod\'elisation, LPSM,\\ 4, Place Jussieu, 75252 Paris cedex 05, France\\
  \printead{e3}
}
\end{aug}

\begin{abstract}
This document is a book chapter which gives a partial survey on post hoc approaches to false positive control.
\end{abstract}

\maketitle

\tableofcontents
\medskip

\input{BNR-chap_body}

\bibliographystyle{abbrv}
\bibliography{BNR-chap}

\end{document}

%% file: BNR-chap_preamble.tex
\usepackage[utf8]{inputenc}
\usepackage{mathtools}
\newtheorem{theorem}{Theorem}[chapter]
\newtheorem{exercise}{Exercise}[chapter]

\newtheorem{remark}{Remark}[chapter]
\newtheorem{lemma}{Lemma}[chapter]

\newcommand{\R}{\mathbb{R}}
\newcommand{\VBonf}{V^{\mbox{\tiny $k_0$Bonf}}}
\newcommand{\VSimes}{V^{\mbox{\tiny Sim}}}
\newcommand{\Rfam}{\mathfrak{R}}

\newcommand{\telque}{\,:\,}

\renewcommand{\P}{\mathbb{P}}

\newcommand{\mtc}{\mathcal}
\newcommand{\wt}[1]{{\widetilde{#1}}}
\newcommand{\wh}[1]{{\widehat{#1}}}
\newcommand{\ol}[1]{\overline{#1}}
\newcommand{\ind}[1]{{\mathbf{1}\left\{#1\right\}}}
\newcommand{\paren}[1]{\left(#1\right)}

\newcommand{\set}[1]{\left\{#1\right\}}
\newcommand{\abs}[1]{\left| #1 \right|}
\newcommand{\cA}{{\mtc{A}}}
\newcommand{\cH}{{\mtc{H}}}
\newcommand{\cX}{{\mtc{X}}}
\newcommand{\Nm}{\mathbb{N}_m}
\newcommand{\NN}{\mathbb{N}}
\newcommand{\Pro}[1]{\mathbb{P}\left(#1\right)}

\graphicspath{
  {./}
  {figures/}
  {chapters/BNR-chap/figures/}
}


%% file: BNR-chap_body.tex
\chapter{On agnostic post hoc approaches to false positive control}

Classical approaches to multiple testing grant control over the amount
of false positives for a specific method prescribing the set of
rejected hypotheses. On the other hand, in practice many users tend to
deviate from a strictly prescribed multiple testing method and follow
ad-hoc rejection rules, tune some parameters by hand, compare several
methods and pick from their results the one that suits them best, etc.
This will invalidate standard statistical guarantees because of the
selection effect. To compensate for any form of such "data snooping",
an approach which has garnered significant interest recently is to derive "user-agnostic", or post hoc, bounds on the false positives
valid uniformly over all possible rejection sets; this allows
arbitrary data snooping from the user.
In this chapter,
we start from a common approach to
post hoc bounds taking into account the $p$-value level sets
for any candidate rejection set, and explain how to calibrate the
bound under different assumption concerning the distribution of $p$-values.
We then build towards 
a general approach to this problem using 
a family of candidate rejection
subsets (call this a reference family) together with associated bounds on the
number of false positives they contain, the latter holding uniformly over the
family. It is then possible to
interpolate from this reference family to find a bound valid for
any candidate rejection subset.
This general program
encompasses in particular the $p$-value
level sets considered initially in the chapter; we illustrate its interest 
in a different context 
where the reference subsets are fixed and
spatially structured. 
These methods are then applied to a genomic example of differential expression study. 
In this chapter, all references are gathered in Section~\ref{notes}.

\section{Setting and basic assumptions}

Let us observe a random variable $X$ with distribution $P$ belonging to some model $\mtc{P}$. Consider  $m$ null hypotheses $H_{0,i} \subset \mtc{P} $, $i \in \Nm = \set{1,\ldots,m}$, for $P$. We denote $\cH_0(P)=\{i \in \Nm \telque P \mbox{ satisfies } H_{0,i}\}$ the set of true null hypotheses and $\cH_1(P)=\Nm\backslash \cH_0(P)$ its complement.
We assume that a $p$-value $p_i(X)$ is available for each null hypothesis $H_{0,i}$, for each $i \in \Nm$. 

We introduce the following assumptions on the distribution $P$, that will be useful in the sequel:
 \begin{align}
&\forall i \in\mathcal{H}_0(P), \:\:\forall t\in[0,1], \:\:\Pro{p_i(X)\leq t} \leq t;
\tag{Superunif}
\label{eq:superunif}\\
&\{p_i(X)\}_{i\in\mathcal{H}_0(P)} \text{ is a family of indep. variables, indep. of } \{p_i(X)\}_{i\in\mathcal{H}_1(P)}.
\tag{Indep}
\label{eq:indep}
\end{align}

\section{From confidence bounds ...}\label{sec:confbound}

Consider some fixed deterministic $S\subset \Nm$. A $(1-\alpha)$-confidence bound $V=V(X)$ for $|S\cap \cH_0(P)| $, the number of false positives  in $S$, is such that
$$
\forall P\in\mtc{P},\, \qquad \P_{X\sim P}\Big(|S\cap \cH_0(P)| \leq {V} \Big)\geq 1-\alpha.
$$

A first example is given by the $k_0$-Bonferroni bound
$
{V}(X)=\sum_{i\in S} \ind{p_i(X)\geq \alpha k_0 /|S|}+k_0-1,
$
for some fixed $k_0\in\Nm$ such that $k_0 \leq |S|$
(otherwise the bound is trivial).
The coverage probability is ensured under \eqref{eq:superunif} by the Markov inequality:
\begin{multline*}
  \P(|S\cap \cH_0(P)| \geq {V}+1)\\
\begin{aligned}
&\leq  \P\bigg(|S\cap \cH_0(P)| \geq \sum_{i\in S\cap \cH_0(P)} \ind{p_i(X)\geq  \alpha k_0/|S|} + k_0\bigg)\\
&=\P\bigg( \sum_{i\in S\cap \cH_0(P)} \ind{p_i(X)< \alpha k_0/|S|} \geq k_0\bigg)\\
&\leq \frac{|S\cap \cH_0(P)|\alpha k_0/|S| }{k_0}  \leq \alpha.
\end{aligned}
\end{multline*}

However, in practice, $S$ is often chosen by the user and possibly depends on the same data set, then denoted $\wh{S}$ to emphasize this dependence;
it typically corresponds to items of potential strong interest. The most archetypal example is when $\wh{S}$ consists of the $s_0$ smallest $p$-values $p_{(1:m)},\dots, p_{(s_0:m)}$, for some fixed value of $s_0 \in\Nm$. In that case, it is easy to check that the above bound does not have the correct coverage: for instance, when the $p$-values are i.i.d. $U(0,1)$ and $\cH_0(P)=\Nm$, we have for $k_0\leq s_0$ (that is, when the bound is informative),
\begin{align*}
\P(|\wh{S}\cap \cH_0(P)| \leq {V})
  &=\P\Big(k_0-1+\textstyle{\sum_{i\in \hat{S}}} \ind{p_i(X)\geq  \alpha k_0/s_0} \geq s_0 \Big)\\
  &=\P\Big(\textstyle{\sum_{i\in \hat{S}}}
  \ind{p_i(X)< \alpha k_0/s_0} \leq k_0-1 \Big)\\
  &= \P\left(p_{(k_0 :m)}(X) \geq \alpha k_0/s_0 \right)\\
  &= \P(\beta(k_0,m-k_0+1)\geq \alpha k_0/s_0),\\
\end{align*}
where $\beta(k_0,m-k_0+1)$ denotes the usual beta distribution with parameters $k_0$ and $m-k_0+1$.
For instance, taking $s_0=10$, $k_0=5$, $\alpha=0.05$ and $m=500$, the latter is approximately equal to $0.005$, while the intended target is $1-\alpha=0.95$.

This phenomenon is often referred to as the selection effect: after some data driven selection, the probabilities change and thus the usual statistical inferences are not valid.

\section{... to post hoc bounds}

To circumvent the selection effect, one way is to aim for a function ${V}(X,\cdot) : S \subset \Nm \mapsto {V}(X,S)\in \mathbb{N}$ (denoted by ${V}(S)$ for short) satisfying %
\begin{equation}\label{aim}
\forall P\in\mtc{P},\, \qquad \P_{X\sim P}\Big(\forall S \subset 
\Nm ,\:
|S\cap \cH_0(P)| \leq {V}(S) \Big)\geq 1-\alpha,
\end{equation}
that is, a $(1-\alpha)$ confidence bound that is valid {\it uniformly} over all subsets $S \subset 
\Nm$.  As a result, for any particular algorithm $\wh{S}$, inequality \eqref{aim} entails $\P\big(
|\hat{S}\cap \cH_0(P)| \leq {V}(\hat{S}) \big)\geq 1-\alpha$, and thus does not suffer from the selection effect. Such a bound
 will be referred to as a {\it $(1-\alpha)$-post hoc confidence bound} throughout this chapter, "post hoc" meaning that 
 the set $S$ can be chosen after having seen the data, and possibly using the data several times.

As a first example, the {\it $k_0$-Bonferroni post hoc bound} is 
\begin{equation}
  \label{eq:k0-bonf}
\VBonf(S)=|S|\wedge\left(\sum_{i\in S} \ind{p_i(X)\geq  \alpha k_0  /m}+k_0-1\right).
\end{equation}
Following the same reasoning as above, it has a coverage at least $1-\alpha$ under \eqref{eq:superunif}:
\begin{multline*}
  \P(\exists S \subset 
\Nm \::\: |S\cap \cH_0(P)| \geq \VBonf(S)+1)\\
\begin{aligned}
&\leq \P\Bigg(\exists S \subset 
\Nm \::\:  \sum_{i\in S\cap \cH_0(P)} \ind{p_i(X)< \alpha k_0 /m} \geq k_0\Bigg)\\
&=\P\Bigg(\sum_{i\in  \cH_0(P)} \ind{p_i(X)< \alpha k_0 /m} \geq k_0\Bigg)\\
&\leq \frac{| \cH_0(P)|\alpha k_0/m }{k_0}  \leq \alpha.
\end{aligned}
\end{multline*}

\begin{remark}
Compared to the $k_0$-Bonferroni confidence bound of Section~\ref{sec:confbound}, $\alpha$ has been replaced by $\alpha |S|/m$,
so that the post hoc bound is much more conservative than a (standard, non uniform, $S$ fixed) confidence bound when $|S|/m$ gets small, which is well expected. This scaling factor is the price paid here to make the inference post hoc. We will see in Sections~\ref{sec:threshold} and~\ref{sec:spatial} that it can be diminished when considering bounds of a different nature.
\end{remark}

\begin{exercise}
For $k_0=1$, when the $p$-values are i.i.d. $U(0,1)$ and $\cH_0(P)=\Nm$, prove that the coverage probability of the $k_0$-Bonferroni post hoc bound is at most $(1-\alpha/m)^m$. Does the Bonferroni post hoc bound provide a sharp coverage in that case?
\end{exercise}

The Bonferroni post hoc bound, while it is valid under no assumption on the dependence structure of the $p$-value family,  may be conservative, in the sense that ${V}(S)$ will be large for many subsets $S$.
For instance, one has $\VBonf(S)=|S|$ (trivial bound) for all the sets $S$ such that 
$S\subset \{i\in\Nm\::\:p_i(X)>\alpha k_0 /m\}$. 

The Bonferroni bound can be further improved under some dependence restriction, with the {\it Simes post hoc bound: }
\begin{equation}\label{equ:Simes}
\VSimes(S)=\min_{1\leq k\leq |S|}\left\{\sum_{i\in S} \ind{p_i(X)\geq  \alpha k /m} + k-1\right\} = \min_{1\leq k\leq |S|} \{V^{\mbox{\tiny $k$Bonf}}(S)\}.
\end{equation}
Its coverage can be computed as follows (using arguments similar as above):
\begin{align}
&\P(\exists S \subset 
\Nm \::\: |S\cap \cH_0(P)| \geq \VSimes(S)+1)\nonumber\\
&\leq \P \Bigg(\exists S \subset 
\Nm \:,\:\exists k\in \{1,\dots,m\} \::\:  \sum_{i\in S\cap \cH_0(P)} \ind{p_i(X)< \alpha k/m} \geq k\Bigg)\nonumber\\
&=\P\left(\exists k\in \{1,\dots,|\cH_0(P)|\} \::\: p_{(k:\cH_0(P))}< \alpha k/m\right)  \label{equ:reasoningSimes}.
\end{align}
Under \eqref{eq:superunif} and \eqref{eq:indep}, this
is lower than or equal to $\alpha|\cH_0(P)|/m\leq \alpha$ by using the {\it Simes inequality}.
More generally, the Simes post-hoc bound is valid in any setting where the Simes inequality holds. This is the case under a specific positive dependence assumption called Positive Regression Dependency on a Subset of hypotheses (PRDS), which is also the assumption under which the Benjamini-Hochberg (BH) procedure has been shown to control the false discovery rate (FDR).

While it uses more stringent assumptions, $\VSimes(S)$ can be much less conservative than $\VBonf$. For instance, if $S=\{i\in\Nm\::\:5\alpha /m \leq p_i(X)< 10 \alpha /m\}$, we have ${V}^{\mbox{\tiny $5$Bonf}}(S)=|S|$ and $\VSimes(S)\leq |S|\wedge 9$, which can lead to a substantial improvement.

\begin{figure}[htb]
\centerline{\includegraphics[scale=0.4]{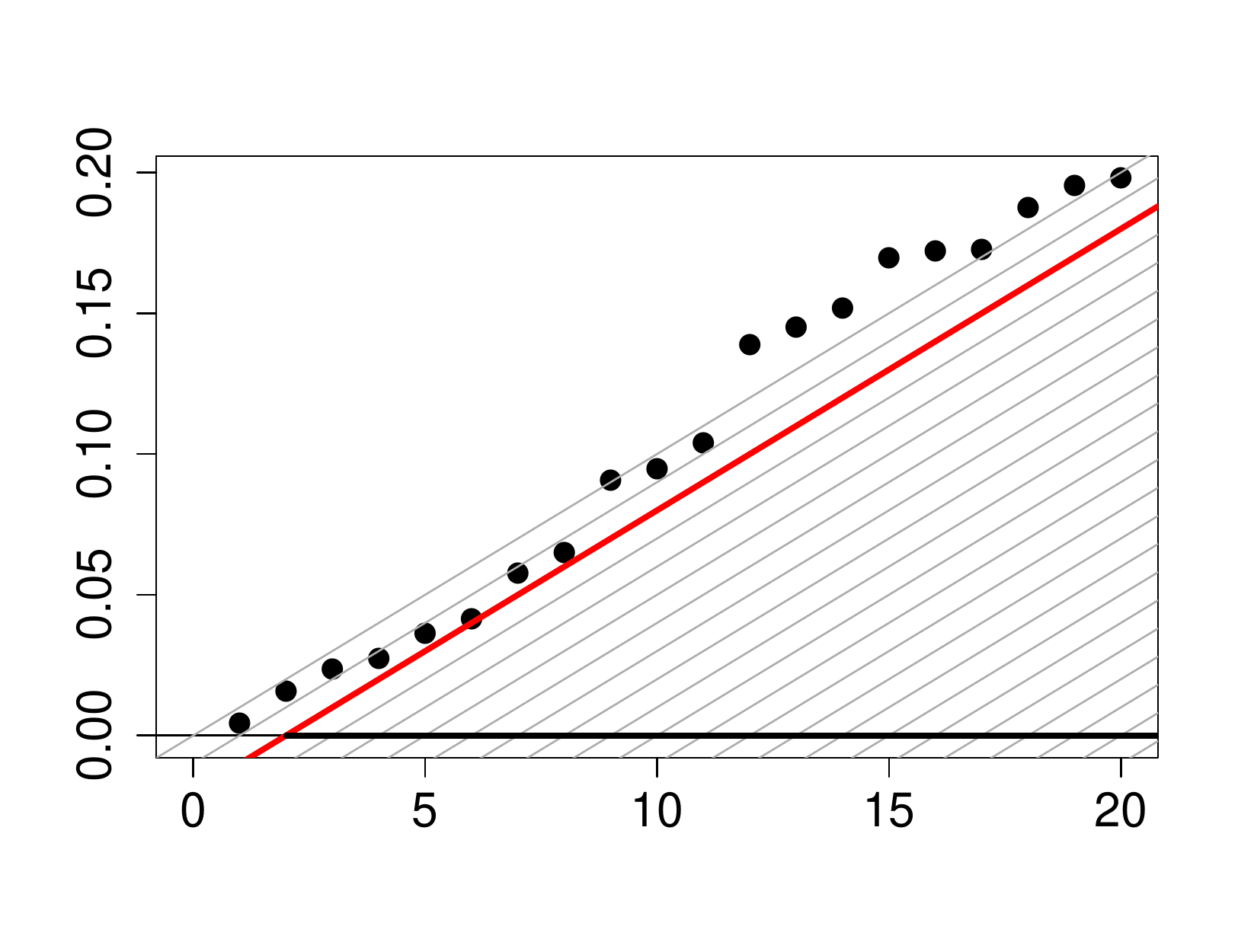}\includegraphics[scale=0.4]{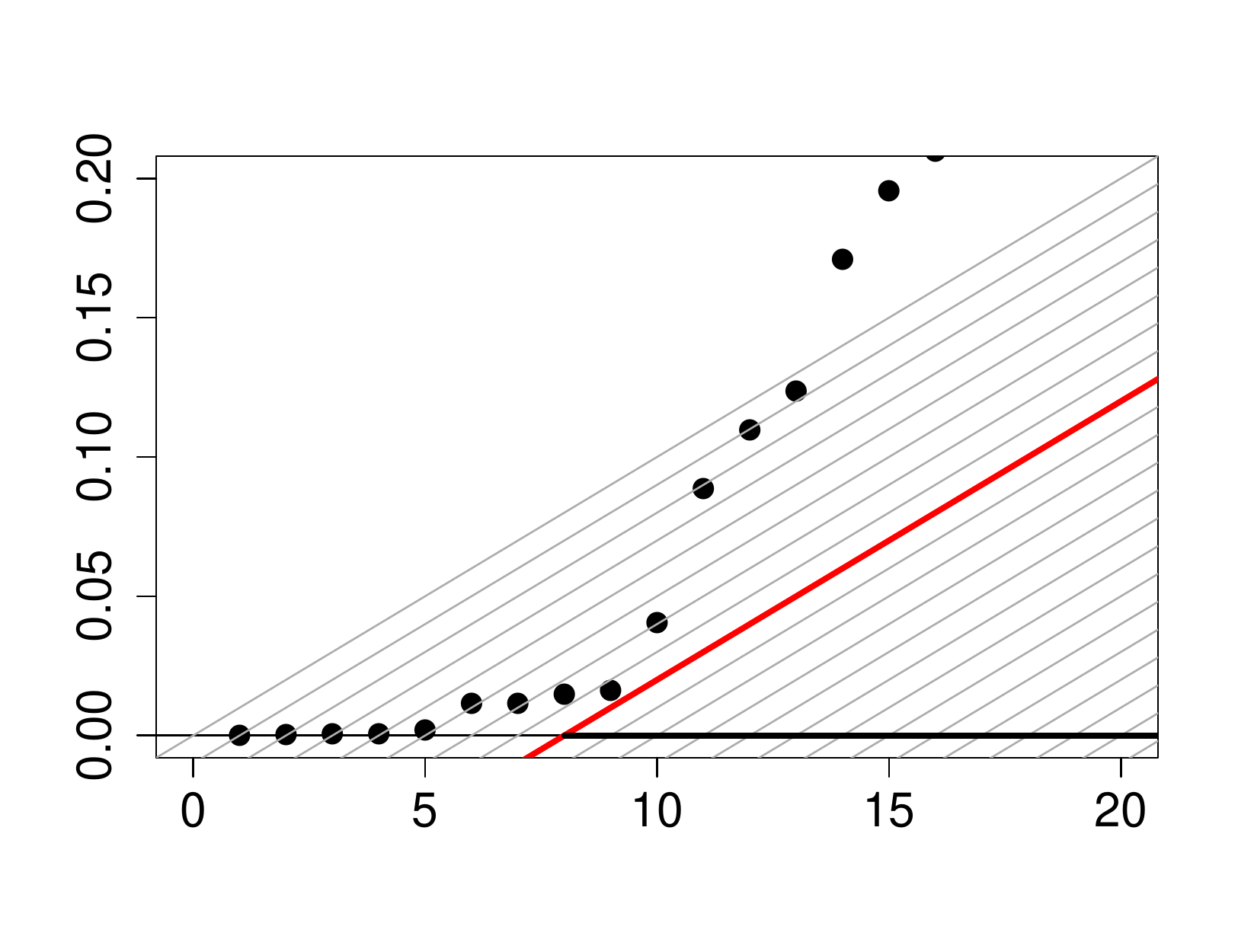}}
\caption[Illustrations for Simes post hoc bound]{Illustration of the Simes post hoc bound \eqref{equ:Simes} according to the expression \eqref{equ:Simeseasy}, for two subsets of $ \Nm$ (left display/right display), both of cardinal $20$ 
and  for $m=50$. The level is $\alpha=0.5$ (taken large only for illustration purposes). Black dots: sorted $p$-values in the respective subsets. Lines: thresholds $k\in\{u+1,\dots,|S|\} \mapsto \alpha(k-u)/m$ (in red for $u = |S|-V(S)$, in light gray otherwise).   
  The post hoc bound $\VSimes(S)$ corresponds the length of the bold line on the $X$-axis.
\label{fig:Simes}}
\end{figure}

From Exercise~\ref{exo:Simeseasy} below, the Simes bound has a nice graphical interpretation:  $|S|-\VSimes(S)$ can be interpreted as the smallest integer $u$ for which the shifted line $v\mapsto \alpha (v-u)/m$ is strictly below the ordered $p$-value curve, see Figure~\ref{fig:Simes}.

\begin{exercise}\label{exo:Simeseasy}
Prove that for all $S\subset \Nm$, $ |S|-\VSimes(S)$ is equal to
 \begin{align}\label{equ:Simeseasy}
\min\{ u \in \{0,\dots,|S|\}\::\:  \forall v\in \{ u+1,\dots,|S|\} \::\: p_{(v:S)} \geq \alpha (v-u)/m\},
\end{align}
where $p_{(1:S)},\dots,p_{(|S|:S)}$ denote the ordered $p$-values of $\{p_i(X),i\in S\}$. [Hint: start from $|S|-\VSimes(S) \leq u$ for some $u$, and find
  equivalent expressions.]
\end{exercise}

\begin{exercise}
  In Figure~\ref{fig:Simes}, check that  $\VSimes(S)=18$ (resp. $\VSimes(S)=12$) in the left (resp. right) situation. Compare to $\VBonf$ for $k_0=7$ ($m=50$, $\alpha=0.5$). 
\end{exercise}

The Simes post hoc bound \eqref{equ:Simes}  has, however, several limitations: first, the coverage is only valid when the Simes inequality holds. This imposes restrictive conditions on the model used, which are rarely met or provable in practice. Second, even in that case, the bound does not incorporate the dependence structure, which may yield conservativeness (see Exercise~\ref{exo:Simes} below). Finally, this bound intrinsically compares the ordered $p$-values to the threshold $k\mapsto \alpha k/m$ (possibly shifted). We can legitimately ask whether taking a different threshold (called template below) does not provide a better bound.

\begin{figure}[htb]
\centerline{\includegraphics[scale=0.7]{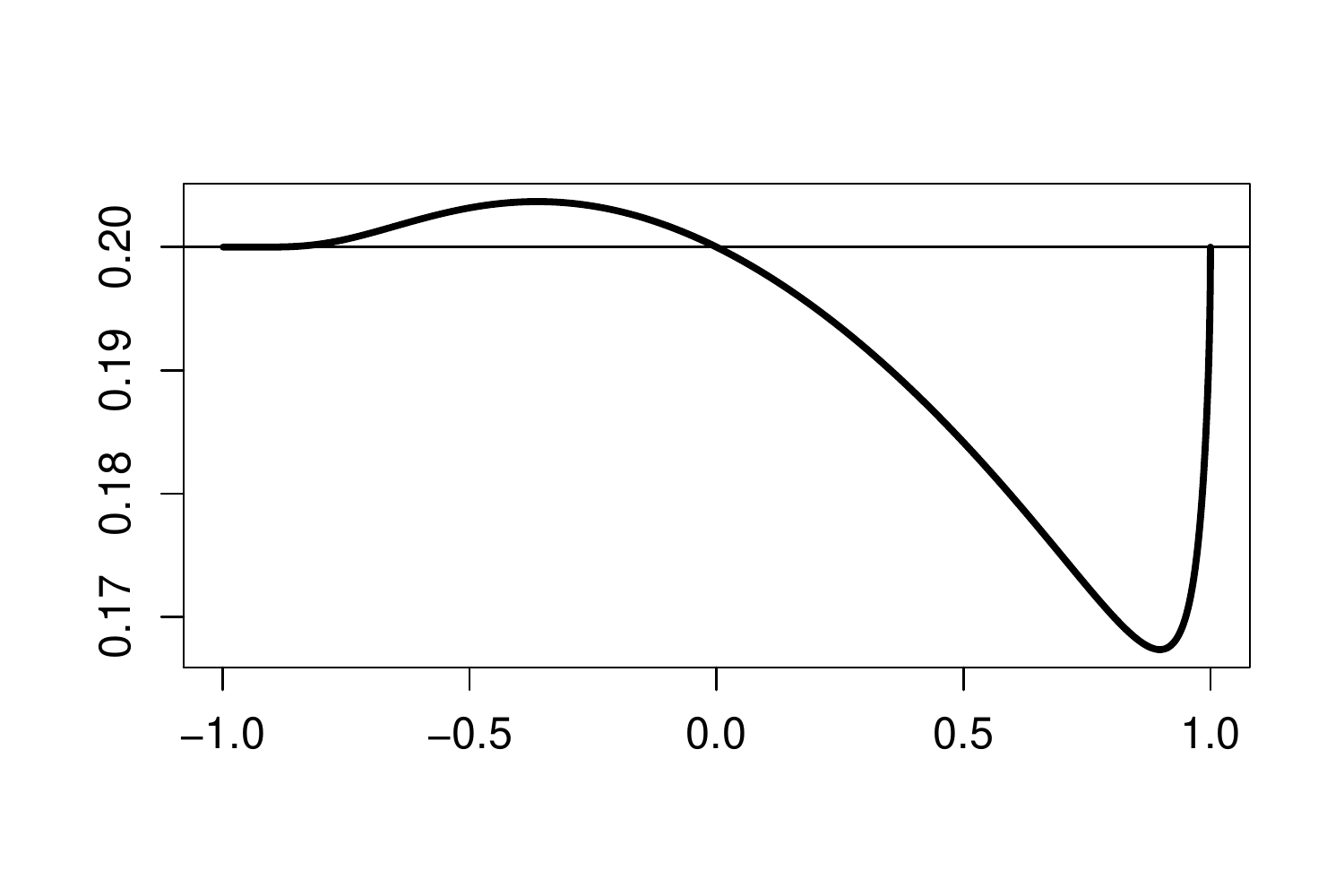}}
\caption{Coverage of the Simes post hoc bound \eqref{equ:Psi} in the setting of Exercise~\ref{exo:Simes} as a function of $\rho$ and for $\alpha=0.2$. \label{fig:simesconserv}}
\end{figure}

\begin{exercise}\label{exo:Simes}
Consider the case $\cH_0(P)=\Nm$, for which $m$ is even, and denote $\ol{\Phi}$ the upper-tail distribution function of a standard $\mathcal{N}(0,1)$ variable. Consider the one-sided testing situation where $p_i= \ol{\Phi}(X_1)$, $1\leq i \leq m/2$ and  $p_i= \ol{\Phi}(X_2)$, $m/2+1 \leq i \leq m$, for a $2$-dimensional Gaussian vector $(X_1,X_2)$  that is centered, with covariance matrix having $1$ as diagonal elements and $\rho\in [-1,1]$ as off-diagonal elements. Show that the coverage probability of the Simes post hoc bound is equal to 
\begin{equation}\label{equ:Psi}
\alpha/2 + \int_{\alpha/2}^\alpha \ol{\Phi}\left(\frac{\ol{\Phi}^{-1}(\alpha)-\rho \ol{\Phi}^{-1}(w)}{(1-\rho^2)^{1/2}}\right)dw +  \int_{\alpha}^\infty \ol{\Phi}\left(\frac{\ol{\Phi}^{-1}(\alpha/2)-\rho \ol{\Phi}^{-1}(w)}{(1-\rho^2)^{1/2}}\right)dw
\end{equation}
The above quantity is displayed in Figure~\ref{fig:simesconserv} for $\alpha=0.2$, as a function of $\rho$. 
\end{exercise}

\section{Threshold-based post hoc bounds}\label{sec:threshold}

\begin{figure}[htb]
\vspace{-10mm}
\centerline{\includegraphics[scale=0.7]{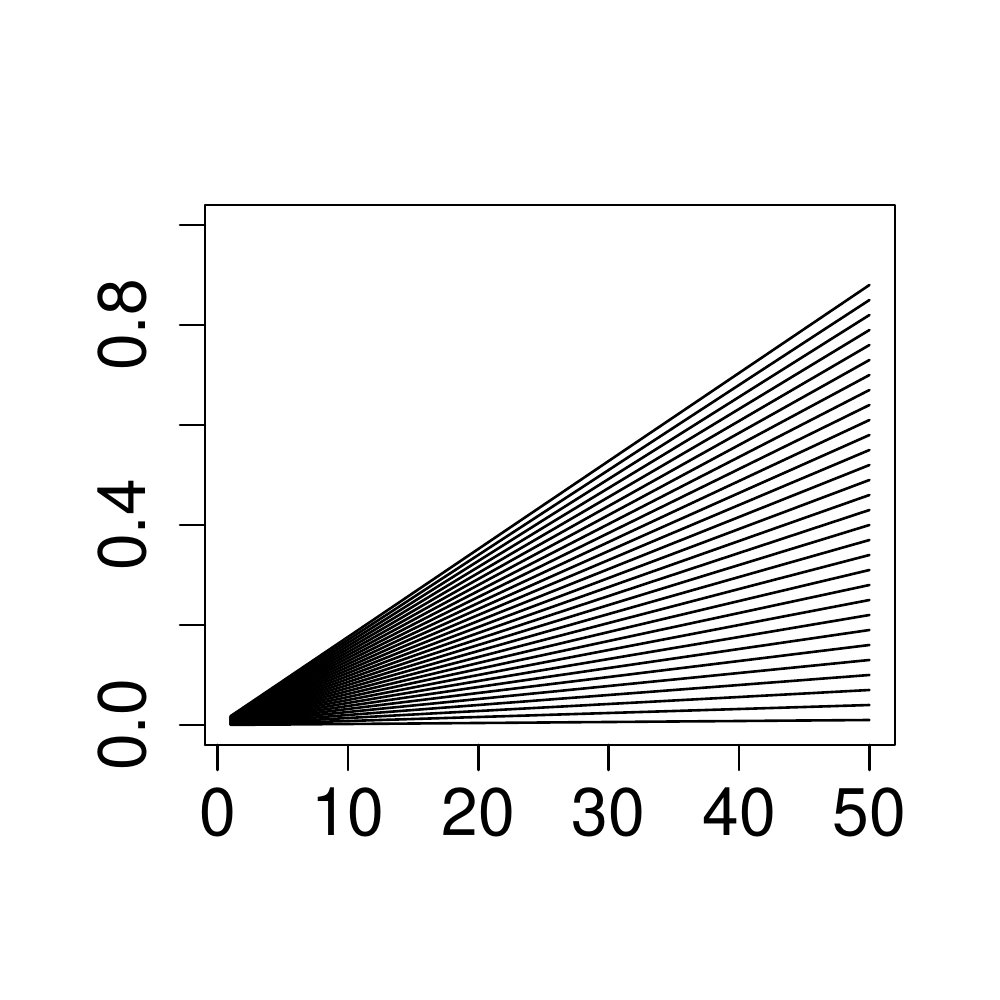}\hspace{-10mm}\includegraphics[scale=0.7]{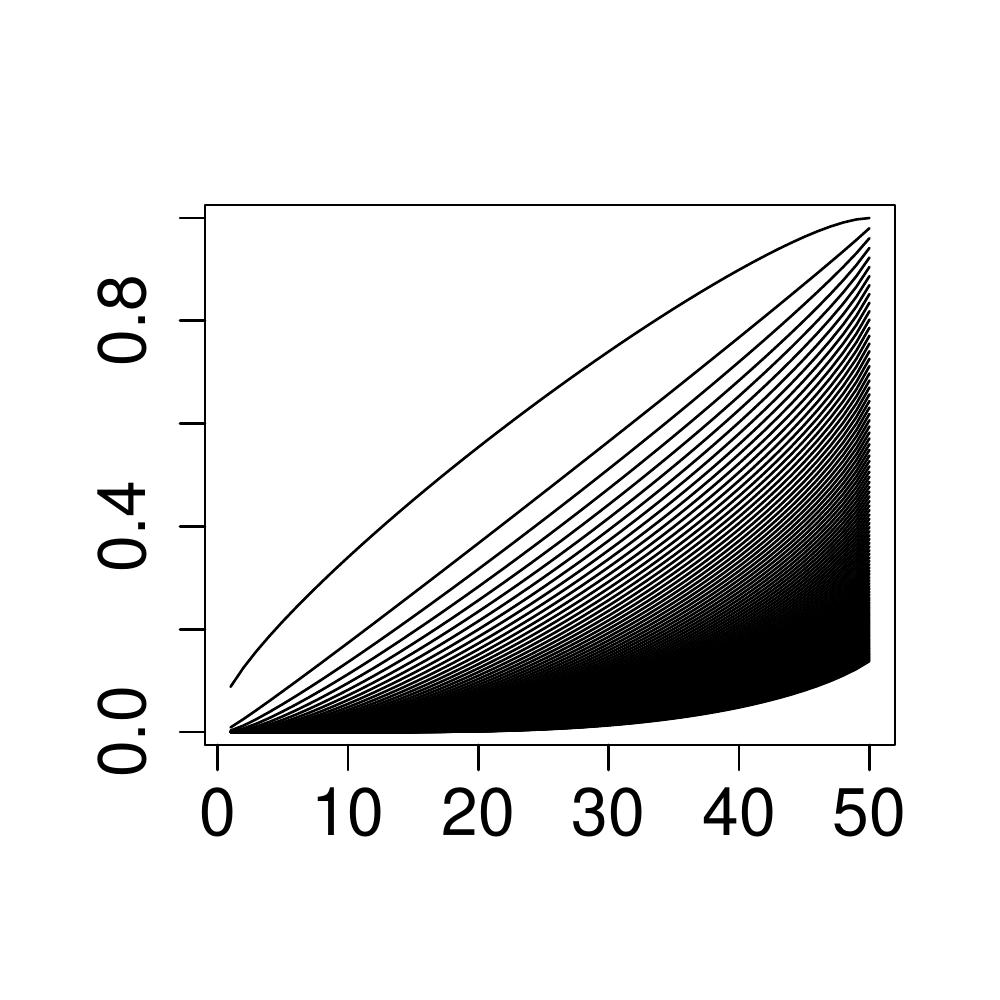}}
\vspace{-7mm}
\caption[Template]{Curves $k\mapsto t_k(\lambda)$ for a wide range of $\lambda$ values. Left: linear template. Right: beta template.
\label{fig:template}}
\end{figure}

More generally, let us consider bounds of the form
\begin{equation}\label{equ:thresholdbounds}
{V}^{\lambda}(S)=\min_{1\leq k\leq |S|} \left\{\sum_{i\in S} \ind{p_i(X)\geq  t_k(\lambda)} + k-1\right\}, \:\: \lambda\in[0,1],
\end{equation}
where  $t_k(\lambda)$, $\lambda\in[0,1]$, $1\leq k \leq m$, is a family of functions, called a template. A template can be seen as a spectrum of curves, parametrized by $\lambda$.
We focus here on the two following examples:
\begin{itemize}
\item Linear template: $t_k(\lambda)=\lambda k/m$, $t_k^{-1}(y)=y m/k$;
\item Beta template: $t_k(\lambda)=$$\lambda$-quantile of $\beta(k,m-k+1)$, $t_k^{-1}(y)=$ $\P(\beta(k,m-k+1)\leq y)$.
\end{itemize}
An illustration for the above templates is provided in Figure~\ref{fig:template}.

For a fixed template, the idea is now to choose one of these curves, that is, one value of the parameter $\lambda=\lambda(\alpha)$, so that the overall coverage is larger than $1-\alpha$. 
Following exactly the same reasoning as the one leading to \eqref{equ:reasoningSimes}, we obtain
\newsavebox{\LHS}\savebox{\LHS}{blahsdbsdgskj sdjkg sdjks sdgjksdjk sdgjk}
\begin{align}
\MoveEqLeft[5]
  {\P(\exists S \subset 
  \Nm \::\: |S\cap \cH_0(P)| \geq {V}^{\lambda}(S)+1)}
  \nonumber\\
    &\leq \P\left(\exists k\in \{1,\dots,|\cH_0(P)|\} \::\: p_{(k:\cH_0(P))}< t_k(\lambda)\right)  \label{equ:reasoningGeneral0}\\
&=\P\left(\min_{k\in \{1,\dots,|\cH_0(P)|\}}\left\{ t_k^{-1}(p_{(k:\cH_0(P))})\right\} < \lambda\right)
 \label{equ:reasoningGeneral},
\end{align}
by letting $t_k^{-1}(y)=\max\{x\in[ 0,1] \::\: t_k(x)\leq y\}$ the generalized inverse of $t_k$ (in general, this is valid provided that for all $k\in\{1,\dots,m\}$, $t_k(0)=0$ and $t_k(\cdot)$ is non-decreasing and left-continuous on $[0,1]$, as in the case of the two above examples).
What remains to be done is thus to calibrate $\lambda=\lambda(\alpha, X)$ such that the quantity \eqref{equ:reasoningGeneral} is below $\alpha$. 

Several approaches can be used for this. It is possible that for the model
under consideration, the joint distribution of $(p_{i}(X))_{i \in \cH_0(P)}$ is
equal to the restriction of some known, fixed distribution on $[0,1]^{\Nm}$ to the coordinates of $\cH_0(P)$ (this is a version of the so-called subset-pivotality condition).
It is met under condition \eqref{eq:indep},
but it is also possible
that the dependence structure of the $p$-values is known (for example, in genome-wide association studies, the structure and strength of linkage disequilibrium can be tabulated from previous studies and give rise to a precise dependence model). In such a situation, the calibration of $\lambda=\lambda(\alpha, X)$ can be obtained either by exact computation, numerical approximation or Monte-Carlo approximation under the full null.

Another situation of interest, on which we focus for the remainder of
this section, is when the null corresponds to an invariant distribution with respect to a certain group of data transformations, which is
the setting for (generalized) permutation tests, allowing for the use of an exact randomization technique. More precisely, assume the existence of a finite transformation group $\mathcal{G}$ acting onto the observation set $\mathcal{X}$.
By denoting $p_{ \cH_0}(x)$ the  null $p$-value vector $(p_i(x))_{i\in \cH_0(P)}$ for  $x\in\mathcal{X}$, we assume that the joint distribution of the transformed null $p$-values is invariant under the action of any $g\in \mathcal{G}$, that is,
\begin{equation}
\forall P \in \mathcal{P},\:\: \forall g \in\mathcal{G},\:\: (p_{\cH_0}(g'.X))_{g'\in \mathcal{G}} \sim (p_{\cH_0}(g'.g.X))_{g'\in \mathcal{G}}  \label{rand} \tag{Rand},
\end{equation}
where $g.X$ denotes $X$ that has been transformed by $g$. 

Let us consider a (random) $B-$tuple $(g_1,g_2,\dots,g_B)$ of $\mathcal{G}$ (for some $B\geq 2$), where $g_1$ is the identity element of $\mathcal{G}$ and $g_2,\dots,g_B$ have been drawn  (independently of the other variables) as i.i.d. variables, each being uniformly distributed on $\mathcal{G}$.
Now, let for all $x\in \mathcal{X}$, $
\Psi(x)= \min_{1\leq k \leq m}\left\{ t_k^{-1}\left(p_{(k:m)}(x)\right)\right\}
$
and consider $\lambda(\alpha, X)=\Psi_{(\lfloor \alpha B\rfloor+1)}$  where $\Psi_{(1)}\leq \Psi_{(2)}\leq \cdots \leq \Psi_{(B)}$ denote the ordered sample $(\Psi(g_j.X), 1\leq j \leq B)$. The following result holds.

\begin{theorem}
\label{thm:lambda-cal}
Under \eqref{rand}, for any deterministic template, the bound ${V}^{\lambda(\alpha,X)}$ is a post hoc bound of coverage $1-\alpha$.
This level is to be understood as a joint
probability with respect to the data and the draw of the group elements $(g_i)_{2\leq i \leq B}$.
\end{theorem}

As a case in point, let us consider a two-sample framework where 
$$X=(X^{(1)},\dots, X^{(n_1)},X^{(n_1+1)},\dots, X^{(n_1+n_2)})\in (\R^m)^{n}$$ is composed of $n=n_1+n_2$ independent $m$-dimensional real random vectors with $X^{(j)}$, $1\leq j \leq n_1$, i.i.d.  $\mathcal{N}(\theta^{(1)},\Sigma)$ (case) and $X^{(j)}$, $n_1+1\leq j \leq n$, i.i.d. $\mathcal{N}(\theta^{(2)},\Sigma)$ (control). Then we aim at testing the null hypotheses $H_{0,i}:$ ``$\theta^{(1)}_i=\theta^{(2)}_i$", simultaneously for $1\leq i \leq m$, without knowing the covariance matrix $\Sigma$. Consider any family of  $p$-values $(p_i(X))_{1 \leq i \leq m}$ such that $p_i(X)$ only depends on the $i$-th coordinate
 $(X_i^{(j)})_{1\leq j \leq n}$ of the observations (e.g., based on difference of the
coordinate means of the two groups).
 Note that $p_{\cH_0}(X)$ is thus a measurable function of $(X_i^{(j)})_{i\in \cH_0,1\leq j \leq n}$. Now, the group $\mtc{G}$ of permutations  of $\{1,\dots,n\}$ is naturally acting on $\cX=(\R^m)^{n}$ via the permutation of the individuals: for all $\sigma \in\mtc{G}$,
$$
\sigma.X= (X^{(\sigma(1))},\dots, X^{(\sigma(n_1))},X^{(\sigma(n_1+1))},\dots, X^{(\sigma(n))}).
$$

\begin{exercise}
Show that $
(X_i^{(1)})_{i\in \cH_0},\dots, (X_i^{(n)})_{i\in \cH_0} 
$ are i.i.d. and prove \eqref{rand}. 
\end{exercise}

\begin{figure}[htb]
\vspace{-5mm}
\centerline{\includegraphics[scale=0.45]{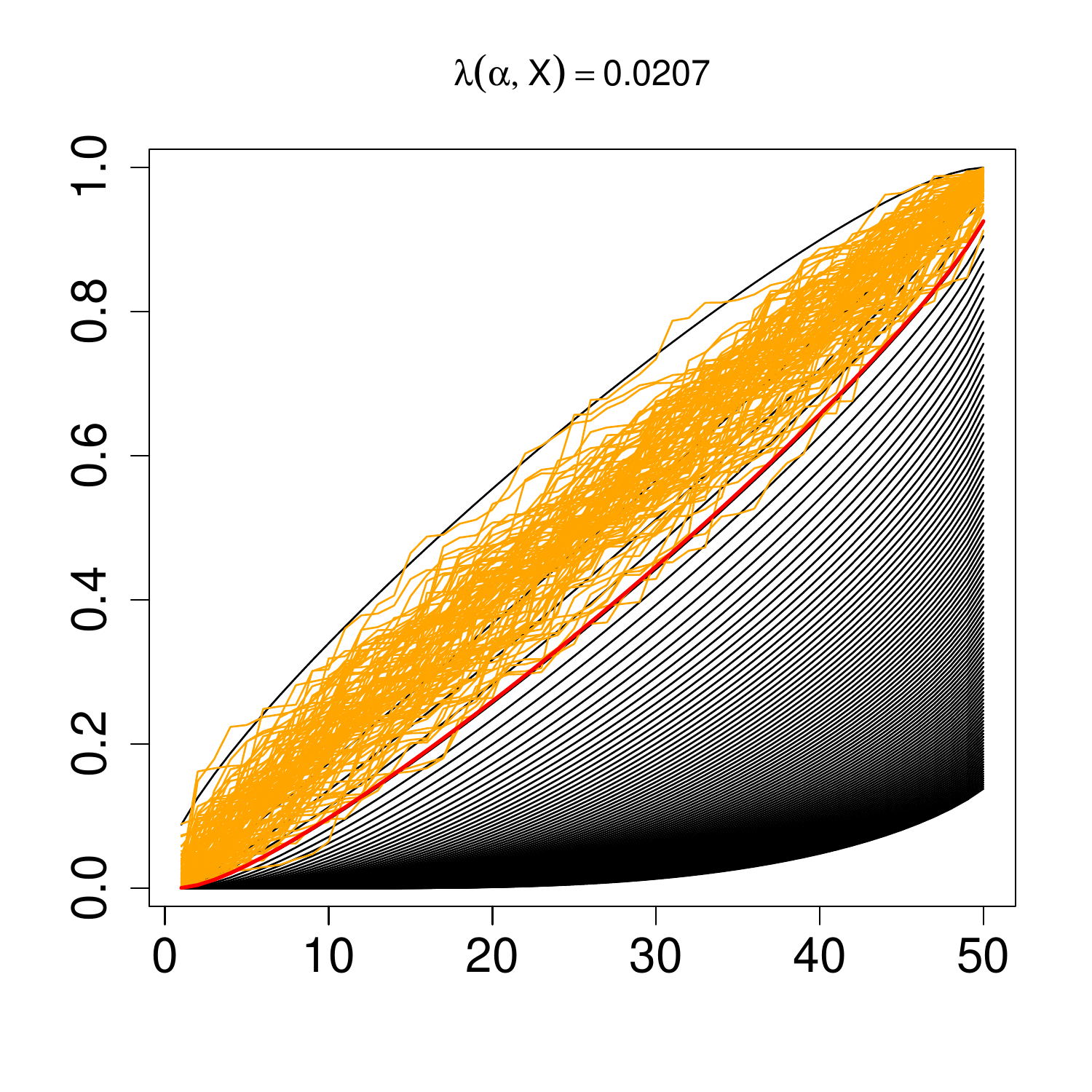}\hspace{-0.6cm}\includegraphics[scale=0.45]{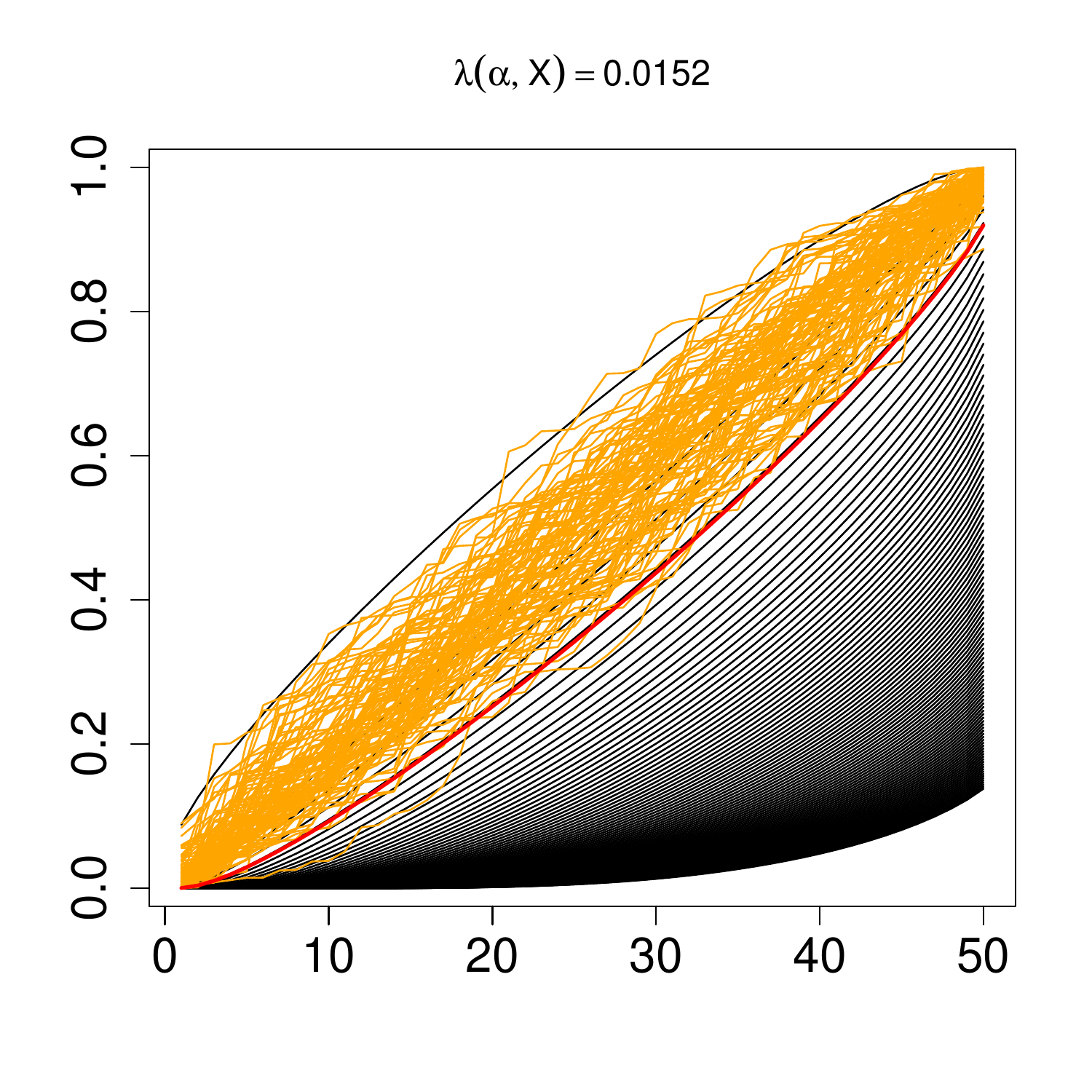}}
\vspace{-7mm}
\caption[$\lambda$-calibration]{Illustration of the $\lambda=\lambda(\alpha,X)$ calibration method on one realization of the data $X$. Black curves: beta template $k\mapsto t_k(\lambda)$ for some range of $\lambda$  values. Orange curves: ordered $p$-values (after permutation) $k\mapsto p_{(k:m)}(g_j.X)$ for $1\leq j \leq B=1000$.  Red curve: $k\mapsto t_k(\lambda(\alpha,X))$. Left panel : full null, right panel : half of true nulls (see text). (Parameters  $m=50$, $\alpha=0.2$, $n_1=50$, $n_2=50$, $\delta=3$.) 
\label{fig:lambdacal}}
\end{figure}

 An illustration of the above $\lambda$-calibration method is provided in Figure~\ref{fig:lambdacal} in the  case where  $\Sigma=I_m$, 
 $$p_i(X)=2\left(1-\Phi\Bigg(s_{n_1,n_2}^{-1}\bigg|n_2^{-1}\sum_{j=n_1+1}^{n_1+n_2} X_i^{(j)}-n_1^{-1}\sum_{j=1}^{n_1} X_i^{(j)}\bigg|\Bigg)\right), $$
 for $s_{n_1,n_2}=(n_1^{-1}+n_2^{-1})^{1/2}$
and using a beta template. 
In the left panel (full null), we have $\theta^{(1)}=\theta^{(2)}=0$, so that $\cH_0(P)=\Nm$. In the right panel (half of true nulls), we have $\theta_i^{(1)}=\theta_i^{(2)}=0$ for $1\leq i\leq m/2$ and $\theta_i^{(1)}=0$, $\theta_i^{(2)}=\delta/s_{n_1,n_2}$ for $m/2+1\leq i\leq m$, for some $\delta>0$, so that $\cH_0(P)=\{1,\dots,m/2\}$. Following expression \eqref{equ:reasoningGeneral0}, $k\mapsto t_k(\lambda(\alpha,X))$ is the highest beta curve such that at most $B \alpha$ orange curves have a point situated below it. This also shows that the above $\lambda$-calibration is slightly more severe when part of the data follows the alternative distribution. This is a commonly observed phenomenon: although the permutation approach
  is valid even when part of the null hypotheses are false, their inclusion in the
  permutation procedure tends to yield test statistics that exhibit more variation
  under permutation, thus inducing more conservativeness in the calibration.

\section{Reference families}\label{sec:R}

We cast the previous bounds in a more general setting, where $(1-\alpha)$--post hoc bounds are explicitly based on a {\it reference family} with some {\it joint error rate} (JER in short) controlling property. This general point of view offers more flexibility and allows us to consider post hoc bounds of a different nature, as for instance those incorporating a spatial structure, see Section~\ref{sec:spatial}.

In general, a reference family is defined by a collection $\Rfam=\big((R_1(X),\zeta_1(X)),$ $\ldots,$ $(R_K(X),\zeta_K(X))\big)$, where the $R_k$'s are data-dependent subsets  of $\Nm$ and the $\zeta_k$'s are data dependent integer numbers (we will often omit the dependence in $X$ to ease notation). The reference  family $\Rfam$  is said to control the JER at level $\alpha$ if 
\begin{equation} 
\forall P\in\mtc{P},\, \qquad\P_{X\sim P}(\forall k \in \NN_K\::\; |R_k(X)\cap \cH_0| \leq \zeta_k(X))\geq 1-\alpha.
\label{eq:jr}
\end{equation}
Markedly, \eqref{eq:jr} is similar to \eqref{aim}, but restricted to some subsets $R_k$, $k\in\NN_K$. The rationale behind this approach is that, while the choice of $S$ is let completely free in \eqref{aim} (to accommodate any choice of the practitioner), the choice of the $R_k$'s and $\zeta_k$'s in \eqref{eq:jr} is done by the statistician and is part of the procedure. Once we obtain a reference family $\Rfam$ satisfying \eqref{eq:jr}, we obtain a post hoc bound by interpolation:
\begin{equation}\label{eq:optbound1}
  V^*_\Rfam(S) = \max\{ \abs{S\cap A}\::\:  A \subset \Nm, \forall k\in \NN_K, |R_k\cap A| \leq \zeta_k \}, \:\:\: S\subset \Nm\,.
\end{equation}
We call $V^*_\Rfam$ the {\it optimal post hoc bound} (built upon the reference family $\Rfam$). 
Computing the bound $V^*_\Rfam(S)$ can be time-consuming,
it actually has NP-hard complexity in a general configuration.
We can introduce the following computable relaxations: for $S\subset \Nm$,
\begin{align}
\ol{V}_\Rfam(S) &=\min_{k \in \NN_K} \paren{ \abs{S\setminus R_k} + \zeta_k} \wedge |S| ;
\label{eq:aug}\\
\wt{V}_\Rfam(S) &=
 \Bigg(\sum_{k \in \NN_K}  \abs{S \cap R_k}\wedge  \zeta_k  + \bigg|S \setminus \bigcup_{k\in \NN_K} R_k\bigg|\Bigg)\wedge |S| . \label{eq:disj}
\end{align}

\begin{exercise}
Show that $V^*_\Rfam(S)\leq \ol{V}_\Rfam(S)$ and $V^*_\Rfam(S)\leq \wt{V}_\Rfam(S)$ for all $S\subset \Nm$.
Moreover, provided that \eqref{eq:jr} holds, show that $V^*_\Rfam$, $\ol{V}_\Rfam$ and $\wt{V}_\Rfam$ are all valid $(1-\alpha)$--post hoc bounds.
\end{exercise}

In addition, the following result shows that the relaxed versions coincide with the optimal bound if the reference sets have some special structure:
\begin{lemma}~
  
\begin{itemize}
\item In the nested case, that is, $R_k\subset R_{k+1}$, for $1\leq k\leq K-1$, we have $\ol{V}_\Rfam = V^*_\Rfam$;
\item In the disjoint case,  that is, $R_k\cap R_{k'}=\emptyset$ for $1\leq k\neq k'\leq K$, we have $\wt{V}_\Rfam = V^*_\Rfam$.
\end{itemize}
\end{lemma}

We can briefly revisit the post-hoc bounds of the previous sections in this general
framework. The $k_0$-Bonferroni post hoc bound~\eqref{eq:k0-bonf} derives from
the one-element reference family $(R=\set{i \in \Nm: p_i(X) <  \alpha k_0/m},\zeta=k_0-1)$. The Simes post hoc bound~\eqref{equ:Simes} derives from
the reference family comprising the latter reference sets for all $k_0 \in \Nm$.
More generally, the threshold-based post hoc bounds $V^\lambda$ of the form \eqref{equ:thresholdbounds} are equal to the optimal bound $V^*_\Rfam$ with $R_k=\{i\in\Nm\::\:p_i(X)< t_k(\lambda)\}$ and $\zeta_k=k-1$, $k \in \Nm$(indeed, these reference sets are nested, so that $ V^*_\Rfam=\ol{V}_\Rfam $). 

How to choose a suitable reference family in general? A general rule of thumb is to choose the reference sets $R_k$ of the same qualitative form as the sets $S$ for which the bound is expected to be accurate. For instance, the Simes post hoc bound will be more accurate for sets $S$ with the smallest $p$-values. In Section~\ref{sec:spatial}, we will choose reference sets $R_k$ with a spatial structure, which will produce a post hoc bound more tailored for spatially structured subsets $S$. 

\section{Case of a fixed single reference set}\label{sec:single}

It is useful to focus first on the case of a single fixed (non-random) reference set $R_1$, with (random) $\zeta_1$ satisfying \eqref{eq:jr}, that is, 
$$\P( | \cH_0(P)\cap R_1| \leq \zeta_1(X))\geq 1-\alpha.$$
(In contrast with the $k_0$-Bonferroni bound~\eqref{eq:k0-bonf} where $\zeta$ was fixed
and $R$ variable, here $R_1$ is fixed and $\zeta_1$ is variable.)
In other words, $\zeta_1(X)$ is a $(1-\alpha)$--confidence bound of $| \cH_0(P)\cap R_1|$.
Several example of such $\zeta_1(X)$ can be build, under various assumptions.

\begin{exercise}\label{exo:confbound}
For $R_1\subset \Nm$ fixed, show that the following bounds are $(1-\alpha)$--confidence bounds for $| \cH_0(P)\cap R_1|$:
\begin{itemize}
\item under \eqref{eq:superunif},  for some fixed $t\in(0,\alpha)$, 
\begin{equation}\label{Markovbounds}
\zeta_1(X)=|R_1|\wedge\Bigg\lfloor \sum_{i\in R_1} \ind{p_i(X)>t}/(1-t/\alpha) \Bigg\rfloor,
\end{equation}
where $\lfloor x\rfloor$ denotes the largest integer smaller than or equal to $x$.
[Hint: use the Markov inequality.]
\item under \eqref{eq:superunif} and \eqref{eq:indep},
\begin{equation}\label{DKWbounds}
\zeta_1(X)=|R_1|\wedge \min_{t\in[0,1)} \left\lfloor \frac{C}{2(1-t)} +
\left(\frac{C^2}{4(1-t)^2} + \frac{\sum_{i\in R_1} \mathbf{1}{\{p_i(X) > t\}}}{1-t}\right)^{1/2} \right\rfloor^2,     
\end{equation}
where  $C=\sqrt{\frac{1}{2}\log \left(\frac{1}{\alpha} \right) }$.
[Hint: use the DKW inequality, that is, for any integer $n\geq 1$, for $U_1,\dots,U_n$ i.i.d. $U(0,1)$, we have $n^{-1} \sum_{i=1}^n \mathbf{1}{\{U_i > t\}} - (1-t) \geq -\sqrt{\log (1/\lambda)/(2n)}$  for all $t\in[0,1]$ with probability at least $1-\lambda$.]
\end{itemize}
\end{exercise}

In addition to the two above bounds \eqref{Markovbounds} and \eqref{DKWbounds}, we can elaborate another bound in the generalized permutation testing framework
~\eqref{rand}, as described in Section~\ref{sec:threshold}. Applying the result of that section, the following bound is also valid:
\begin{equation}\label{Betabounds}
\zeta_1(X)=\min_{1\leq k\leq |R_1|} \left\{\sum_{i\in R_1} \ind{p_i(X)\geq  t_k(\lambda(\alpha,X))} + k-1\right\},
\end{equation}
where $t_k(\lambda)$ denotes the $\lambda$-quantile of a $\beta(k,|R_1|-k+1)$ distribution and $\lambda(\alpha,X)=\Psi_{(\lfloor \alpha B\rfloor+1)}$, where $\Psi_{(1)}\leq \Psi_{(2)}\leq \cdots \leq \Psi_{(B)}$ denote the ordered sample $(\Psi(g_j.X), 1\leq j \leq B)$ for which $\Psi(x)= \min_{1\leq k \leq |R_1|}\left\{ t_k^{-1}\left(p_{(k:|R_1|)}(x)\right)\right\}$ (see the $\lambda$-calibration method of Section~\ref{sec:threshold}).\\

Once a proper choice of $\zeta_1(X)$ has been done,  the optimal post hoc bound can be computed as follows: for any $S\subset \Nm$, 
$
V^*_\Rfam(S)=\ol{V}_\Rfam(S) = \wt{V}_\Rfam(S)=|S\cap R_1^c|+ \zeta_1(X)\wedge |S\cap R_1|.
$
When $S$ is large and does not contain very small $p$-values, this bound can be sharper than the Simes bound.

\begin{exercise}\label{exo:single}
Let us consider a post bound based on the single reference family $R_1=\Nm$ and $\zeta_1(X)$ as in \eqref{DKWbounds} (choosing $t=1/2$).
For $S$ such that $S\subset  \{i\in \Nm\::\: p_i(X)> \alpha |S|/m\}$, show that $\VSimes(S)=|S|$ and $V^*_\Rfam(S)=|S|\wedge\zeta_1(X)\leq |S|\wedge 2\left(\log \left(\frac{1}{\alpha} \right) + 2\sum_{i\in \Nm} \ind{p_i(X)>  1/2} \right)$.
\end{exercise}

Finally, while the case of a single reference set can be considered as an elementary example, the  bounds developed in this section will be useful in the next section, for which several fixed reference sets $R_k$ are considered, and thus several (random) $\zeta_k$ should be designed.

\section{Case of spatially structured reference sets}\label{sec:spatial}

We consider here the case where the null hypotheses $H_{0,i}$, $1\leq i\leq m$, have a spatial structure, 
and we are interested in obtaining accurate bounds on $|S\cap \cH_0(P)|$ for subsets $S$ of the form $S=\{i\in\Nm\::\:i_0\leq i\leq j_0\}$, for some $1\leq i_0<j_0\leq m$. 

In that case, it is natural to choose $R_k$ formed of contiguous indices.
To be concrete, consider
reference sets consisting of disjoint intervals of the same size : assume $m=Ks$ for some integers $K>0$ and $s>0$ and let
\begin{equation}
  \label{eq:ref-spatial}
R_k=\{(k-1)s +1,\dots,ks\},\: k\in \NN_K.
\end{equation}
When each of these regions is considered in isolation, Section~\ref{sec:single} suggested several approaches (in the appropriate settings \eqref{eq:superunif}, \eqref{eq:indep} or \eqref{rand}) of a specific form $\zeta_k(X)=f(R_k,\alpha,X)$, to underline the dependence of $\zeta_k(X)$ in $R_k$ and $\alpha$.
By using a simple union bound, it is then straightforward to show that the JER control \eqref{eq:jr} is satisfied for 
\begin{equation}\label{deviceunionbound}
\zeta_k(X)=f(R_k,\alpha/K,X),\:  k\in \NN_K.
\end{equation}
When the reference regions $R_k$ are disjoint as in the example~\eqref{eq:ref-spatial}
above, we can use the proxy $\wt{V}_\Rfam(S)$ (see~\eqref{eq:disj}) which is known
to coincide with the optimal bound $ V^*_\Rfam(S)$.
This gives rise to a post hoc bound  that accounts for the spatial structure of the data.
 
 \begin{exercise}\label{exo:influenzeofK}
Compute $\zeta_k(X)$ in the case where $\zeta_1(X)=f(R_1,\alpha,X)$ is given by \eqref{Markovbounds} ($t=\alpha^2$) and \eqref{DKWbounds}. In each case, for a given $k$, discuss how $\zeta_k(X)$ fluctuates when the size of the family $K$ increases.
 \end{exercise}

 When considering the reference regions defined by segments~\eqref{eq:ref-spatial},
we have to prescribe a scale ($s$ here, the size of the segments). It is possible 
to extend this to a multi-scale approach, choosing overlapping reference intervals $R_k$
at different resolutions arranged in a tree structure, where parent sets are formed by taking union of (disjoint) children sets taken at a finer resolution. Furthermore, the proxy~\eqref{eq:disj} has to be replaced
by a more elaborate one, minimizing over all possible multi-scale partitions made of such
reference regions. This can still be computed efficiently by exploiting the
the tree structure.
Doing so, the post hoc bound will be more scale adaptive to sets $S$ with possibly various sizes. The price to pay lies in the cardinality $K$ of the family, which gets larger. However, this does not necessarily make the corresponding bound much larger, as Exercise~\ref{exo:influenzeofK} shows when using the bound \eqref{DKWbounds},
since the level $\alpha$ only enters it logarithmically.


\section{Applications}
\label{sec:appli}

Differential gene expression studies in cancerology aim at identifying genes whose mean expression level differs significantly between two (or more) cancer populations, based on a sample of gene expression measurements from individuals from these populations. We consider here a  microarray data set\footnote{Taken from Chiaretti \emph{et. al.}, {\em Clinical cancer research}, 11(20):7209--7219, 2005.} consisting of expression measurements for more than $12,000$ genes for biological samples from $n = 79$ individuals with B-cell acute lymphoblastic leukemia (ALL). A subset of cardinal $n_1=37$ of these individuals harbor a specific mutation called  BCR/ABL, while the remaining $n_2=42$ don't. One of the goals of this study is to identify those genes for which there is a difference in the mean expression level between the mutated and non-mutated population. This question can be addressed, after relevant data preprocessing, by performing a statistical test of equality in means for each gene. A classical approach is then to derive a list of ``differentially expressed'' genes (DEG) as those passing a FDR correction by the Benjamini-Hochberg (BH) procedure at a user-defined level. For example, 163 genes are selected by the BH procedure at level $\alpha = 0.05$. We note that although the usage of the BH procedure is standard for multiple two-sample tests and widely accepted in the biomedical literature, we have no formal guarantee that it is mathematically justified -- in particular, genes are not independent, and there is no proof that the PRDS assumption holds in this setting.

In this section, we illustrate how the post hoc inference framework introduced in the preceding sections can be applied to this case to build confidence envelopes for the proportion of false positives (Section~\ref{sec:confenv}), and to obtain bounds on data-driven sets of hypotheses (Section~\ref{sec:volcano}), and on sets of hypotheses defined by an a priori structure (Section~\ref{sec:structured}). These numerical results were obtained using the R package \textbf{sansSouci}, version 0.8.1\footnote{Available from \texttt{https://github.com/pneuvial/sanssouci}.}.

\subsection{Confidence envelopes}\label{sec:confenv}

\begin{figure}[!htp]
  \centering
 \includegraphics[height=5cm]{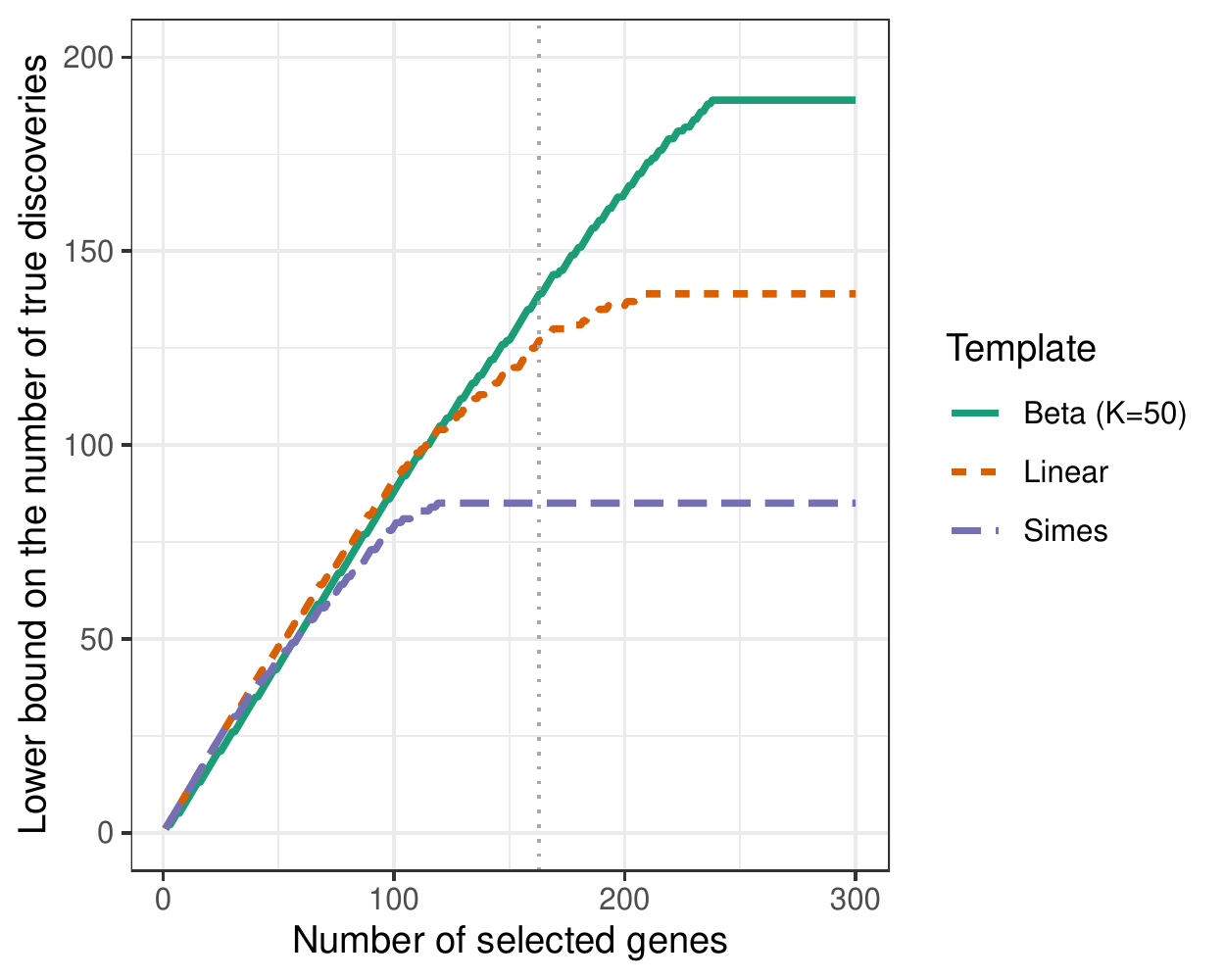}
 \includegraphics[trim=0 0 90 0, clip, height=5cm]{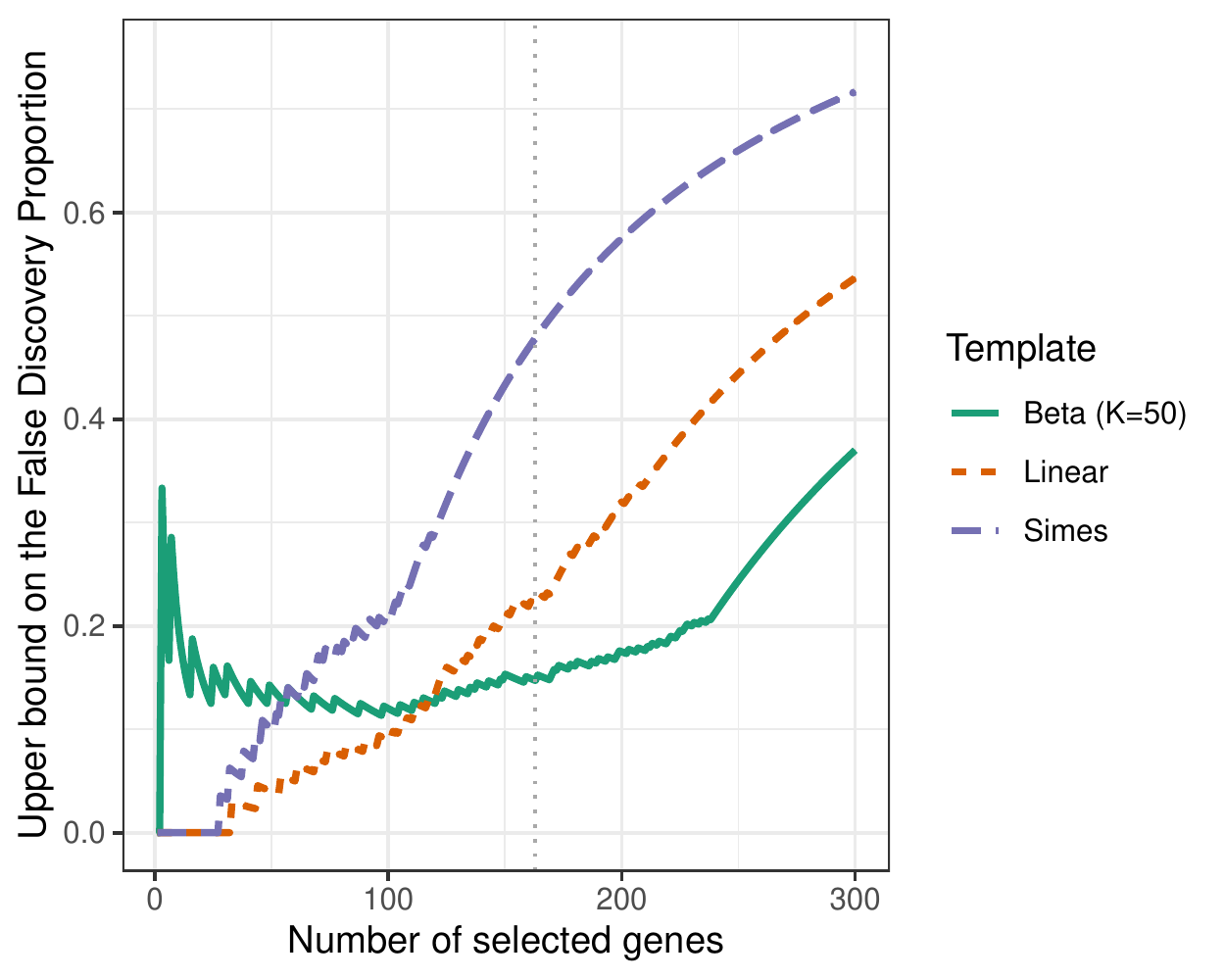}
  \caption{Confidence bounds on the number of true positives (left) and on the proportion of false positives (right) for several reference families: Simes reference family  (long-dashed purple curve), linear template after $\lambda$-calibration (dashed red curve), and beta template after $\lambda$-calibration (solid green curve).}
  \label{fig:confidence-bounds}
\end{figure}

In absence of specific prior information on relevant subsets of hypotheses to consider, it is natural to focus on subsets consisting of the most
significant hypotheses. 
Specifically, we define the $k-$th $p$-value level set $S_k$ as the set of the $k$ most significant hypotheses, corresponding to the $p$-values $(p_{(1:m)}, p_{(2:m)}, \dots , p_{(k:m)})$, and consider post hoc bounds associated to $S_k$ for $k \in \Nm$. Figure~\ref{fig:confidence-bounds} provides \emph{post hoc confidence envelopes}
for the ALL data set, for $\alpha = 0.1$. While $(1-\alpha)$-lower confidence bounds on the number of true positives of the form $\set{\paren{k, \abs{S_k} - \overline{V}(S_k)}: k \in \mathbb{N}_m}$ are displayed in the left panel, $(1-\alpha)$-upper confidence bounds on the proportion of false positives $\set{\paren{k, \overline{V}(S_k)/\abs{S_k}}: k \in \mathbb{N}_m}$ are shown in the right panel.

The confidence envelopes are built from the Simes bound~\eqref{equ:Simes} (long-dashed purple curve), and from two bounds obtained from Theorem~\ref{thm:lambda-cal} by $\lambda$-calibration using $B=1,000$ permutation of the sample labels, based on the two templates introduced in Section~\ref{sec:threshold}: the dashed red curve corresponds to the linear template with $K=m$, and the solid green curve to the beta template with $K=50$. Note that Assumption \eqref{rand} holds because we are in the two-sample framework described after Theorem~\ref{thm:lambda-cal}.
The vertical line in Figure~\ref{fig:confidence-bounds} corresponds to the 163 genes selected by the BH procedure at level $5\%$. The Simes bound ensures that the FDP of this subset is not larger than 0.48. As noted above concerning the BH procedure, we have a priori no guarantee that this bound is valid, because such multiple two-sample testing situations have not been shown to satisfy the PRDS assumption under which the Simes inequality is valid\footnote{In this particular case, $\lambda$-calibration with the linear template yields $\lambda(\alpha) > \alpha$, which a posteriori implies that the Simes inequality was indeed valid.}. In contrast, the $\lambda$-calibrated bounds built by permutation are by construction valid here. Moreover, both are much sharper than the Simes bound while the $\lambda$-calibrated bound using the linear template is twice smaller, ensuring FDP$< 0.23$, and even smaller for the beta template with $K=50$. The bound obtained by $\lambda$-calibration of the linear template is uniformly sharper that the original Simes bound \eqref{equ:Simes}, which corresponds to $\lambda = \alpha$. This illustrates the adaptivity to dependence achieved by $\lambda$-calibration.  The bound obtained from the beta template is less sharp for $p$-value level sets $S_k$ of cardinal less than $k=120$, and then sharper. This is consistent with the shape of the threshold functions displayed in Figure~\ref{fig:template}. 

\subsection{Data-driven sets}\label{sec:volcano}

A common practice in the biomedical literature is to only retain, among the genes called significant after multiple testing correction, those whose ``fold change'' exceeds a prescribed level. The fold change is the ratio between the mean expression levels of the two groups. With the notation of Section \ref{sec:threshold}, the fold-change of gene $i$ is given by $ \Delta_i = \overline{X}_i^{(2)}/\overline{X}_i^{(1)}$, where $\overline{X}_i^{(1)} = n_1^{-1}\sum_{j=1}^{n_1}X_i^{(j)}$ and $\overline{X}_i^{(2)} = n_2^{-1}\sum_{j=1}^{n_2}X_i^{(j)}$.

\begin{figure}[!htp]
  \centering
  \includegraphics[width=8cm]{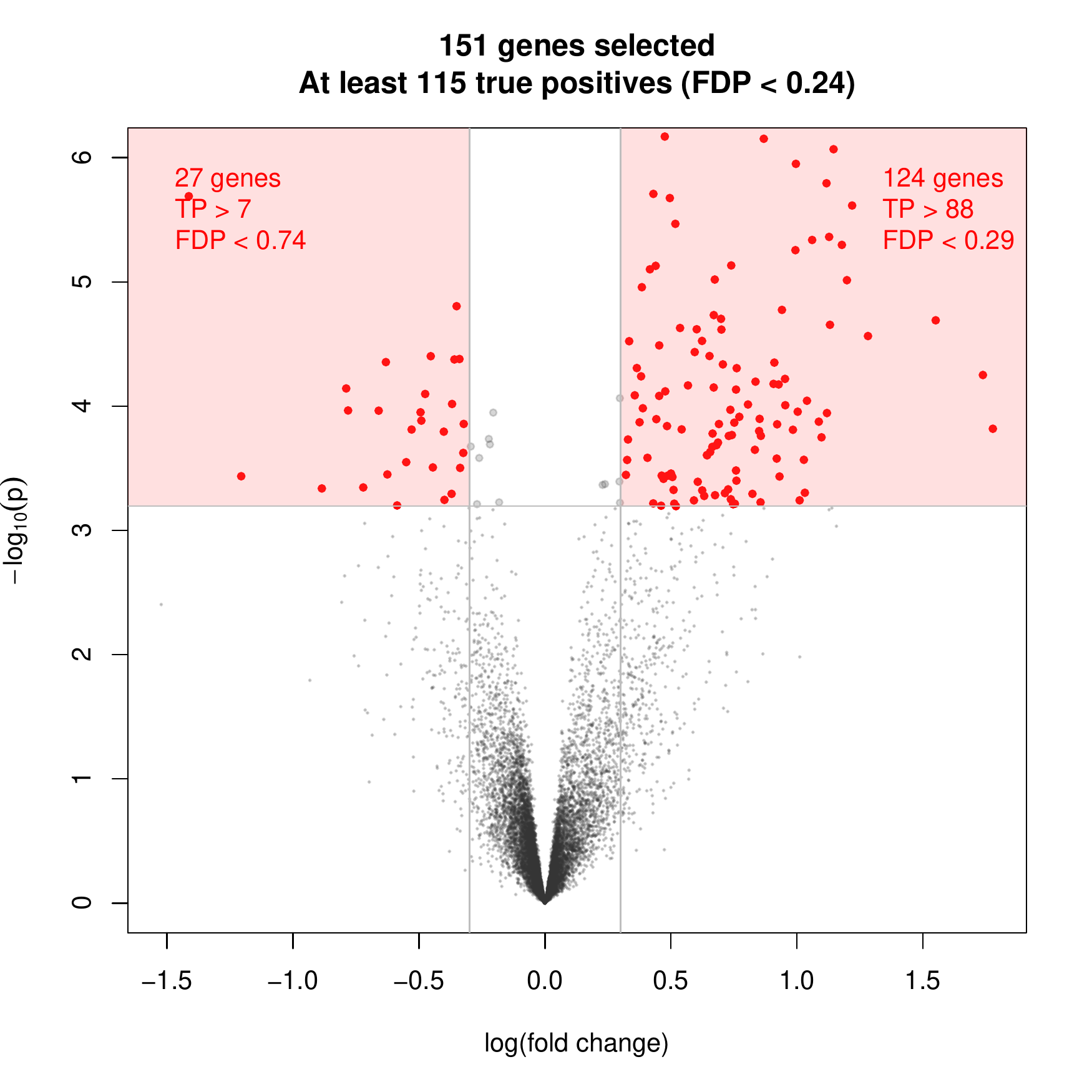}
  \caption{Post-hoc inference for volcano plots}
  \label{fig:volcano-plot}
\end{figure}

 This is illustrated by Figure~\ref{fig:volcano-plot}, where each gene is represented as a point in the ($\log$(fold change), $-\log(p)$) plan. This representation is called a ``volcano plot'' in the biomedical literature.  Among the 163 genes selected by the BH procedure at level 0.05, 151 have an absolute log fold change larger than 0.3. As FDR is not preserved by selection, FDR controlling procedures provide no statistical guarantee on such data-driven lists of hypotheses. 
In contrast, the post hoc bounds proposed in this chapter are valid for such data-driven sets. The two shaded boxes in Figure~\ref{fig:volcano-plot} correspond to the data-driven subsets  $S^{\rm BH} \cap S^-$ and $S^{\rm BH} \cap S^+$, where $S^{\rm BH}$ is the set of 163 genes selected by the BH procedure at level $0.05$,  $S^-= \set{i \in \Nm, \log(\Delta_i) < -0.3}$ and $S^+= \set{i \in \Nm, \log(\Delta_i) > +0.3}$.
The post hoc bounds on the number of true positives in $S^{\rm BH} \cap S^+, S^{\rm BH} \cap S^-$ and  $S^{\rm BH} \cap (S^+ \cup S^-)$ obtained by the Simes bound and by the $\lambda$-calibrated linear and beta templates are given in Table~\ref{tab:volcano}. Both $\lambda$-calibrated bounds are more informative than the Simes bound, in the sense that they provide a higher bound on the number of true confidence. Moreover, they have proven $(1-\alpha)$-coverage, whereas the coverage of the Simes bound is a priori unknown for multiple two-sample tests.  None of the two $\lambda$-calibrated bounds dominates the other one, which is in line with the fact that the linear template is well-adapted to situations with smaller $p$-value level sets than the beta template.

Finally, we also note that the bound obtained for $S^+ \cup S^-$ is systematically larger than the sum of the two individual bounds, which, again, is in accordance with the theory. 
\begin{table}[ht]
\centering
\begin{tabular}{lrrrrrrr}
  \hline
 & n & Simes & Linear & Beta(K=50) \\ 
  \hline
  $S^{\rm BH} \cap S^-$ & 124 & 62 & 88 & 100 \\ 
  $S^{\rm BH} \cap S^+$ & 27 & 1 & 7 & 5 \\ 
  $S^{\rm BH} \cap (S^+ \cup S^-)$ & 151 & 79 & 114 & 127 \\ 
   \hline
\end{tabular} 
\label{tab:volcano}
\caption{Post hoc bounds on the number of true positives in $S^{\rm BH} \cap S^+, S^{\rm BH} \cap S^-$ and  $S^{\rm BH} \cap (S^+ \cup S^-)$ obtained by the post hoc bounds displayed in Figure~\ref{fig:confidence-bounds}.}
\end{table}

\subsection{Structured reference sets}\label{sec:structured}

In this section we give an example of application of the bounds mentioned in Section~\ref{sec:spatial}. Our biological motivation is the fact that gene expression activity can be clustered along the genome. 

The $m$ individual hypotheses are naturally partitioned into $23$ subsets, each corresponding to a given chromosome. Within each chromosome, we consider sets of $s=10$ successive genes as in (\ref{eq:ref-spatial}). Hence, we focus on a reference family with the following elements 
$$
R_{c,k}=\{(k-1)s +1,\dots, \min(ks, m_c)\}, \quad k \in \NN_{K_c}, \quad c \in\{1,\dots,23\},
$$
where, in chromosome $c$, $m_c$ denotes the number of genes, $K_c=\lceil m_c/s \rceil$ the number of corresponding regions. 
In addition, for each $(c,k)$ we use 
$\zeta_{c,k}(X)=f(R_{c,k},\alpha_c/K_c,X)$
coming from the union bound \eqref{deviceunionbound} in combination with the device \eqref{DKWbounds} and $\alpha_c = \alpha m_c/m$. This choice accounts for a union bound over all the chromosomes. 
As shown in Exercise~\ref{exo:confbound}, $\zeta_{c,k}(X)$ is a valid upper confidence bound for $| \cH_0(P)\cap R_{c,k}|$ under \eqref{eq:superunif} and \eqref{eq:indep}. 
In this genomic example, \eqref{eq:indep} may not hold, so we have in fact no formal guarantee that this bound is valid. Therefore, the results obtained below are merely illustrative of the approach and may not have biological relevance.

\begin{figure}[!htp]
  \centering
 \includegraphics[width=\textwidth]{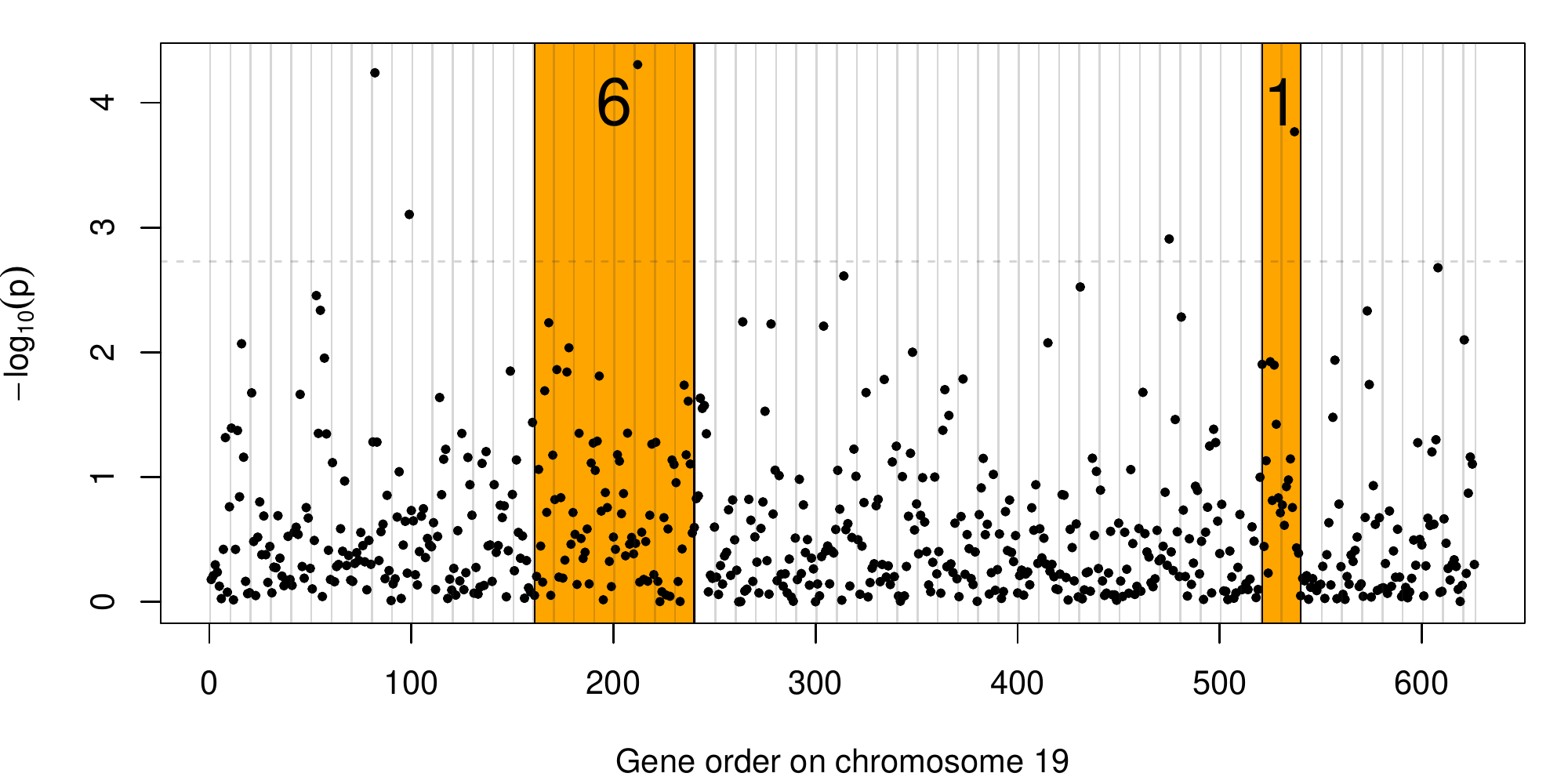}
  \caption{Evidence of locally-structured signal on chromosome 19 detected by the bound \eqref{DKWbounds}.}
  \label{fig:locally-structured}
\end{figure}

We report the results for chromosome $c=19$, which contains $m_c = 626$ genes. In this particular case, we obtain trivial bounds $\zeta_{c,k}(X) = \abs{R_{c,k}}$ for all $k \in \NN_{K_c}$. Therefore, the proxy  $\tilde{V}^*_\mathfrak{R}$ defined in \eqref{eq:disj} for disjoint sets does not identify any signal for this chromosome. However, non-trivial bounds can be obtained via the multi-scale approach briefly mentioned in Section~\ref{sec:spatial}. The idea is to enrich the reference family by recursive binary aggregation of the neighboring $R_{c,k}$. The total number of elements in this family is less than $2K_c$. In our example, it turns out that \eqref{DKWbounds} yields 6 true discoveries in the interval $R_{17:24}$ and 1 true discovery in the interval $R_{53:54}$, where we have denoted
$$R_{u:v}=\bigcup_{u \leq k \leq v} R_{c,k}.$$
This is illustrated by Figure~\ref{fig:locally-structured} where the individual $p$-values are displayed (on the $-\log_{10}$ scale) as a function of their order on chromosome 19. The sets $R_{17:24}$ and  $R_{53:54}$ are highlighted in orange, with the corresponding number of true discoveries marked in each region. 
We obtain a non-trivial bound not because of the large effect of any individual gene, but because of the presence of sufficiently many moderate effects. In particular, in the rightmost orange region in Figure~\ref{fig:locally-structured}, the distribution of $-\log_{10}(p)$ is shifted away from $0$ when compared to the rest of chromosome  19. In comparison, we obtain trivial bounds $\overline{V}_\mathfrak{R}(R_{53:54}) = |R_{53:54}|=2s$ and $\overline{V}_\mathfrak{R}(R_{17:24}) = |R_{17:24}|= 8s$ from (\ref{eq:aug}) both for the linear or the beta template. These numerical results illustrate the interest of the bounds introduced in Section~\ref{sec:spatial} in situations where one expects the signal to be spatially structured.

\section{Bibliographical notes}
\label{notes}

The material exposed in this chapter is mainly a digested account of the article~\cite{BNR2019}. The seminal work~\cite{GW2006} introduced the idea of
false positive bounds for arbitrary
rejection sets. It started from the idea of building a confidence set
on the set of null hypotheses $\cH_0(P)$, and introduced the concepts of
{\em augmentation procedure} and {\em inversion procedure}. The latter consists in
building a confidence set based on
the inversion of tests for $\cH_0(P) = \cA$ for all $\cA \subset \Nm$.
The former starts from a set $R$ with controlled $k$-familywise error rate,
and the proposed associated post hoc bound is~\eqref{eq:jr} (for the one-element reference
family $(R,\zeta=k-1)$). The name {\em augmentation} refers to a similar idea
found in~\cite{DL2008}. The relaxation~\eqref{eq:jr} can in this sense be called
``generalized augmentation procedure''. 
A post hoc
bound for an arbitrary rejection set based on a closed test principle
was proposed in \cite{GS2011}. It can
also be seen as a reformulation of the inversion procedure of~\cite{GW2006}.
Post-hoc bounds over a large class of reference families extracted from classical
FDR control procedures combined with martingale techniques were recently proposed
in~\cite{KatRam2018}.
The principle of the graphical representation used in Figure~\ref{fig:Simes}
to visualize the Simes inequality-based bound originates from J.~Goeman.

The use of generalized permutation procedures in a multiple testing framework
has been explored in several landmark works
\cite{WY1993,RW2005,Mein2006,DL2008,HG2018,HSG2018}. The subset-pivotality condition has been defined in \cite{WY1993}.
Assumption \eqref{rand} has been introduced in \cite{HG2017} and is a weaker version of the randomization hypothesis of \cite{RW2005}.
 The phenomenon of
conservativeness in the permutation-based calibration mentioned at the end of
Section~\ref{sec:threshold}, when not all the null hypotheses are true, can be
in part alleviated by using a step-down principle (see~\cite{RW2005} for a
seminal work on this topic and~\cite{BNR2019} for more details on this approach
in the specific setting considered here). The choice of the size $K$ of the
reference family, which can be crucial in practice, is also discussed
in~\cite{BNR2019}.  

Multiple testing for spatially structured hypotheses is in itself a very active
and broad area of research. 
It has been specifically considered in conjunction with
post-hoc bounds in \cite{MKTG2015}. The use of the reference family approach for
post-hoc bounds in combination with spatially structured hypotheses has been studied
in \cite{DBNR2018}, where the notion of tree- (or forest-)structured reference regions
is introduced, along with an efficient algorithm to compute the optimal bound
$V^*_\Rfam$ in this setting.

The Simes inequality \cite{Sim1986} is a particularly nice and elegant
theoretical device with manifold applications in multiple testing which is still a very active research area, see, e.g., \cite{Block2013,Bodnar2017,Finner2017}.
The DKW inequality with optimal constant was proved in \cite{Mass1990}. 
The Benjamini-Hochberg (BH) procedure has been introduced in \cite{BH1995}, where it is also proved to control  the false discovery rate (FDR). A huge literature on FDR control has followed this seminal paper.

The data used for the application part are taken from \cite{CLGV+2005}. The fact that the signal is clustered along the genome is motivated by previous studies showing possible links between gene expression and DNA copy number changes or other regulation mechanisms~\cite{RSBV+2005,SVRB+2006}.

\section*{Acknowledgements}

This work has been supported by ANR-16-CE40-0019 (SansSouci) and ANR-17-CE40-0001 (BASICS).
The first author acknowledges the support from the german DFG under the Collaborative Research Center SFB-1294 ``Data Assimilation''.
